\documentclass[12pt]
{amsart}

\usepackage{amssymb}
\usepackage{amsmath}
\usepackage{array}

\begingroup\makeatletter\ifx\SetFigFont\undefined
\def\x#1#2#3#4#5#6#7\relax{\def\x{#1#2#3#4#5#6}}%
\expandafter\x\fmtname xxxxxx\relax \def\y{splain}%
\ifx\x\y   
\gdef\SetFigFont#1#2#3{%
  \ifnum #1<17\tiny\else \ifnum #1<20\small\else
  \ifnum #1<24\normalsize\else \ifnum #1<29\large\else
  \ifnum #1<34\Large\else \ifnum #1<41\LARGE\else
     \huge\fi\fi\fi\fi\fi\fi
  \csname #3\endcsname}%
\else
\gdef\SetFigFont#1#2#3{\begingroup
  \count@#1\relax \ifnum 25<\count@\count@25\fi
  \def\x{\endgroup\@setsize\SetFigFont{#2pt}}%
  \expandafter\x
    \csname \romannumeral\the\count@ pt\expandafter\endcsname
    \csname @\romannumeral\the\count@ pt\endcsname
  \csname #3\endcsname}%
\fi
\fi\endgroup

\newcommand{\Civ}[4]{
\begin{picture}(3952,895)(5825,-4136)
\thinlines
\put(6001,-3961){\circle*{336}}
\put(7201,-3961){\circle*{336}}
\put(8401,-3961){\circle*{336}}
\put(9601,-3961){\circle{336}}
\put(6169,-3961){\line( 1, 0){864}}
\put(7369,-3961){\line( 1, 0){864}}
\put(8569,-3911){\line( 1, 0){864}}
\put(8569,-4011){\line( 1, 0){864}}
\put(5941,-3436){\makebox(0,0)[lb]{\smash{\SetFigFont{6}{7.2}{rm}#1}}}
\put(7141,-3436){\makebox(0,0)[lb]{\smash{\SetFigFont{6}{7.2}{rm}#2}}}
\put(8341,-3436){\makebox(0,0)[lb]{\smash{\SetFigFont{6}{7.2}{rm}#3}}}
\put(9541,-3436){\makebox(0,0)[lb]{\smash{\SetFigFont{6}{7.2}{rm}#4}}}
\end{picture}
}

\newcommand{\G}{\mbox{$\mathbf{G}_2$}}
\newcommand{\E}[1]{\mbox{$\mathbf{E}_{#1}$}}
\newcommand{\F}{\mbox{$\mathbf{F}_4$}}

\newcommand{\Pf}{{\em Proof}. }
\newcommand{\EPf}{\hfill$\square$}

\newcommand{\Z}{\mbox{$\mathbf Z$}}
\newcommand{\R}{\mbox{$\mathbf R$}}
\newcommand{\CC}{\mbox{$\mathbf C$}}
\newcommand{\Q}{\mbox{$\mathbf H$}}
\newcommand{\Ca}{\mbox{$\mathbf{Ca}$}}

\newcommand{\SU}[1]{\mbox{$\mathbf{SU}(#1)$}}
\newcommand{\U}[1]{\mbox{$\mathbf{U}(#1)$}}
\newcommand{\SP}[1]{\mbox{$\mathbf{Sp}(#1)$}}
\newcommand{\SO}[1]{\mbox{$\mathbf{SO}(#1)$}}

\newcommand{\Spin}[1]{\mbox{$\mathbf{Spin}(#1)$}}

\newtheorem{thm}{Theorem}

\newtheorem{rmk}{Remark}

\newtheorem{prop}{Proposition}
\newtheorem{lem}{Lemma}
\newtheorem{exs}{Examples}

\title{Taut representations of compact simple Lie groups}
\author{Claudio Gorodski}
\address{Instituto de Matem\'atica e Estat\'\i stica\\Universidade
de S\~ao Paulo\\Rua do Mat\~ao, 1010\\
S\~ao Paulo, SP 05508-090\\
Brasil}
\email{gorodski@ime.usp.br}

\date{\today}

\begin{document}

\begin{abstract}
In this paper we classify the reducible representations of 
compact simple Lie groups all of whose orbits 
are tautly embedded in Euclidean space 
with respect to $\Z_2$ coefficients.
\end{abstract}

\maketitle 

\section{Introduction}

The main result of this paper is the 
following classification theorem.

\begin{thm}\label{thm:taut-red}
A taut reducible representation of a 
compact simple Lie group is one of the 
following representations:
{\small
\[ \begin{array}{|c|c|c|}
\hline
\SU n,\;n\geq3 & \CC^n\oplus\cdots\oplus\CC^n & \mbox{$k$ copies, where $1<k<n$} \\
\hline
\SO n,\;n\geq3, n\neq4 & \R^n\oplus\cdots\oplus\R^n & \mbox{\textrm{$k$ copies, where $1<k$}} \\
\hline
\SP n,\;n\geq1 & \CC^{2n}\oplus\cdots\oplus\CC^{2n} & \mbox{\textrm{$k$ copies, where $1<k$}} \\ 
\hline
\G & \R^7\oplus\R^7 & \mbox{---} \\
\hline
\Spin6 & \R^6\oplus\CC^4 & \R^6=\mathrm{(vector)},\;\CC^4=\mathrm{(spin)} \\
\hline
       & \R^7\oplus\R^8 & \\
\Spin7 & \R^8\oplus\R^8 & \R^7=\mathrm{(vector)},\;\R^8=\mathrm{(spin)} \\
       & \R^8\oplus\R^8\oplus\R^8 & \\
       & \R^7\oplus\R^7\oplus\R^8 & \\
\hline
       & \R^8_0\oplus\R^8_+ & \\
\Spin8 & \R^8_0\oplus\R^8_0\oplus\R^8_+ & \R^8_0=\mathrm{(vector)},\;
\R^8_+=\mathrm{(halfspin)} \\
       & \R^8_0\oplus\R^8_0\oplus\R^8_0\oplus\R^8_+ & \\
\hline
\Spin9 & \R^{16}\oplus\R^{16} & \R^{16}=\mathrm{(spin)} \\ 
\hline
\end{array}\]
}
\end{thm}

A representation of a compact Lie group is called \emph{taut}
if all of its orbits are taut submanifolds of the 
representation space. 
Carter and West introduced in~\cite{CW1} the concept of tautness
for submanifolds (see also~\cite{CR}).
Fix a field of coefficients $F$ (herein assumed to be $\Z_2$). 
Let $M$ be a properly embedded submanifold of 
an Euclidean space~$\R^m$.
For each $p\in\R^m$, consider the squared distance function 
$L_p:M\to\R$ given by $L_p(x)=||x-p||^2$. It is a 
consequence of the Morse index theorem that the critical points
of $L_p$ are nondegenerate, i.e.~$L_p$ is a Morse 
function, if and only if $p$ is not a focal point of $M$.  
Now~$M$ is called \emph{$F$-taut}, or simply $\emph{taut}$, if
$L_p$ is a perfect Morse function for every $p$ in $\R^m$ that is not 
a focal point of $M$. We recall that a Morse function is said to 
be \emph{perfect} if the Morse inequalities are equalities for the 
function restricted to any sublevel set. 
As a consequence
of the proof of the Morse inequalities, one sees that an 
equivalent definition of $F$-tautness for a submanifold
$M\subset \mathbf{R}^m$ is that the induced homomorphism 
\begin{equation*}
 H_*(M\cap B;F)\to H_*(M;F) 
\end{equation*}
in singular homology is injective for almost every closed 
ball $B$ in $\mathbf{R}^m$. It is then clear
that tautness is conformally
invariant.

A compact surface in $\mathbf{R}^3$ which is taut 
is either a round sphere or a cyclide of Dupin (see~\cite{Ba2}); 
the latter can 
all be constructed as the image of a torus of revolution
under a M\"obius transformation. Pinkall and Thorbergsson
found in~\cite{PTh2} the homeomorphism classes of the compact
$3$-dimensional manifolds that admit taut embeddings, 
and their list consists of seven manifolds.
The first three are $S^1\times S^2$ and its quotients
$S^1\times\R P^2$ and $S^1\times_{\mathbf{Z}_2}S^2$. The 
next three are $S^3$ and its quotients $\R P^3$ and
$S^3/\{\pm1,\pm i,\pm j,\pm k\}$ (the so-called 
quaternion space). The last example is the torus $T^3$. 
It follows from
the Chern-Lashof theorem~\cite{CL1} that a taut substantial 
(namely, nor contained in an affine hyperplane)
embedding of a sphere must 
be spherical and of codimension one. 
If $M$ is 
an $n$-dimensional taut hypersurface in $\mathbf{R}^{n+1}$ which has 
the same integral homology as $S^k\times S^{n-k}$, then Cecil and Ryan proved
in~\cite{CR2} that $M$ has precisely two principal curvatures at each point
and that the principal curvatures are constant along the corresponding
curvature distributions.
Bott and Samelson proved in~\cite{BS} that 
the orbits of the isotropy representations of the symmetric spaces,
sometimes called generalized flag manifolds, 
are tautly embedded submanifolds, 
although they did not use this terminology.
The generalized flag manifolds are homogeneous examples of 
submanifolds which belong to another very important, more general
class of submanifolds called isoparametric submanifolds.
Hsiang, Palais and Terng studied in~\cite{HPT}
the topology of isoparametric submanifolds and proved, 
among other things, that they and their focal submanifolds are taut. 

Most of the examples of taut embeddings known are homogeneous spaces. 
In~\cite{Th4} Thorbergsson posed some questions regarding
the problem of which homogeneous spaces admit taut embeddings
and derived some necessary topological conditions for the existence 
of a taut embedding which allowed him to conclude that 
certain homogeneous spaces cannot be tautly embedded
(see also~\cite{He}), among others
the lens spaces distinct from the real projective space.
Many proofs have been given of the tautness
of special cases of generalized flag manifolds where the arguments are 
easier.
No new examples of taut embeddings
of homogeneous spaces besides the generalized flag manifolds 
were known until Gorodski and Thorbergsson classified in~\cite{GTh3}
(see also~\cite{GTh1})
the irreducible representations of compact Lie groups 
all of whose orbits are tautly embedded; we call representations
with this property \emph{taut}. It turns out that
the classification includes three new representations which 
are not isotropy representations of symmetric spaces, thereby 
supplying many new examples of 
tautly embedded homogeneous spaces. 
In~\cite{GTh2} Gorodski and Thorbergsson provided another proof of 
the tautness 
of those orbits by adapting
the proof of Bott and Samelson to that case. 
It is interesting to remark that those three representations 
precisely coincide
with the representations of cohomogeneity three of the compact 
Lie groups which are not orbit equivalent to the isotropy representation
of a symmetric space. (Recall that two representations are said
to be \emph{orbit equivalent} if there is an isometry 
between the representation spaces mapping the orbits of the first
representation onto the orbits of the second one.) 
As mentioned above, in this paper we extend the classification 
in~\cite{GTh3} to the case in which the representation is reducible 
and the group is simple.

The author would like to thank Professor Gudlaugur Thorbergsson 
for useful discussions, and FAPESP and CNPq for 
financial support. 

\section{Preliminary material}

In this section, we collect results that will be used
later to prove that certain representations are or are not taut.
We start with a following simple remark, namely,
\emph{every summand of a taut reducible representation is taut}.
Indeed, this is because an orbit of a summand is also an orbit of the sum,
and it implies that taut reducible representations
are sums of taut irreducible ones. So, in order to classify 
taut reducible representations, we need just
to decide which of those sums are allowed. 
We shall do that for simple groups. 

We begin by recalling the main result of~\cite{GTh3}.

\begin{thm}[\cite{GTh3}]\label{thm:taut-irred}
A taut irreducible
representation of a compact connected Lie group
is either orbit equivalent to the isotropy representation of a symmetric
space or it is one of the following orthogonal representations ($n\geq2$):
{\small
\[ \begin{array}{|c|c|}
\hline
\SO2\times\Spin9 & \mathrm{(vector)}\otimes_{\mathbf R}\mathrm{(spin)} \\
\U2\times\SP n & \mathrm{(vector)}\otimes_{\mathbf C}\mathrm{(vector)} \\
\SU2\times\SP n & \mathrm{(vector)}^3\otimes_{\mathbf H}\mathrm{(vector)} \\
\hline
\end{array}\]
}
\end{thm}

Since the groups appearing in the table of Theorem~\ref{thm:taut-irred}
are nonsimple, now we can refine the remark above and state that
\emph{every summand of a taut reducible representation is 
orbit equivalent to the isotropy representation of a symmetric
space}. Throughout the paper, we shall make use of the tables of 
isotropy representations of a symmetric spaces given in~\cite{Wolf}.
The irreducible representations orbit equivalent to the isotropy 
representation of a symmetric space are also classified (see~\cite{EH1}). 
Lists with some of the principal isotropy subgroups of these representations 
can be found in~\cite{HPT,Straume2}. 

The fundamental result about taut sums of representations is:

\begin{prop}[\cite{GTh3}]\label{prop:taut-sum}
Let $\rho_1$ and $\rho_2$ be representations of a
compact connected Lie group $G$ on
$V_1$ and $V_2$, respectively. Assume that $\rho_1\oplus\rho_2$ is
$F$-taut. Then the restriction of $\rho_2$ to the
isotropy group
$G_{v_1}$ is taut for every $v_1\in V_1$.
Furthermore, we have that $p(G(v_1,v_2);F)=p(Gv_1;F)\,p(G_{v_1}v_2;F)$, where
$p(M;F)$
denotes the Poincar\'e polynomial of $M$ with respect to the field $F$.
In particular, $G_{v_1}v_2$ is connected
and $b_1(G(v_1,v_2);F)=b_1(Gv_1;F)+b_1(G_{v_1}v_2;F)$, 
where $b_1(M;F)$ denotes the first Betti number of $M$ with respect to
$F$.
\end{prop}

We give examples of how Proposition~\ref{prop:taut-sum}
can be used. These are taken from~\cite{GTh1}.

\begin{exs}\label{exs:exs}
\rm (i) Let $G=\SO n$ and let $\rho_1$ be the $\SO n$-conjugation
on the space $V_1$ of real traceless symmetric $n\times n$ matrices. Then
$\rho_1$ is taut since it is the isotropy
representation of the symmetric space ${\bf SL}(n,{\bf R}) /\SO n$. Let
$\rho_2$ be any other nontrivial  representation of
$\SO n$ with representation space $V_2$. Then $\rho_1\oplus\rho_2$ cannot be
taut if $n\ge 3$.  To see this, let $v_1\in V_1$ be
a regular point. Then
$G_{v_1}$ is the discrete group consisting of all diagonal matrices with
determinant one and entries $\pm1$ on the diagonal.
The kernel of $\rho_2$ is contained in the center of $\SO n$. Since $n\ge
3$, we see that $G_{v_1}$ cannot be contained in the
kernel of $\rho_2$. Hence there is an element $v_2\in V_2$ that is not fixed
by $G_{v_1}$. It follows that $G_{v_1}v_2$ is
disconnected. Now Proposition~\ref{prop:taut-sum}
implies that $\rho_1\oplus\rho_2$ is not taut.
The same argument applies more generally whenever 
$\rho_1$ is a taut 
representation of a compact connected 
Lie group $G$ such that its principal isotropy subgroup 
is discrete and not central. 

(ii) Now let $G$ be a compact connected simple Lie group of rank at 
least two and
let $\rho_1$ denote the adjoint representation of
$G$. We assume that $G$ is simply connected. Let
$\rho_2$ be any other nontrivial representation of $G$. Then
$\rho_1\oplus\rho_2$ is not taut. To see this let
$T$ be a maximal torus in $G$. We denote the representation spaces of
$\rho_1$ and $\rho_2$ by $V_1$ and $V_2$ respectively.
There is a regular element $v_1\in V_1$ with $G_{v_1}=T$. The restriction of
$\rho_2$ to $T$ has a discrete kernel that is
contained in the center of $G$. If $v_2\in V_2$ is a $T$-regular point
then the isotropy subgroup $T_{v_2}$ coincides with the
kernel of $\rho_2|T$. Hence $G_{v_1}v_2$ is diffeomorphic to $T$ and it
follows that $b_1(G_{v_1}v_2;F)$ is equal to the rank
of
$G$. In particular $b_1(G_{v_1}v_2;F)\ge 2$.
Now notice that the isotropy group of $(v_1,v_2)$ is also $T_{v_2}$.
Hence $\pi_1(G(v_1,v_2))=T_{v_2}$ which implies $H_1(G(v_1,v_2);{\bf
Z})=T_{v_2}$ since $T_{v_2}$ is Abelian. If $G\ne \Spin {4k}$
then the center of
$G$ is a cyclic group and it follows that
$b_1(G(v_1,v_2);F)\le 1$.
If $G=\Spin{4k}$, then $k\ge2$ and we get $b_1(G_{v_1}v_2;F)=2k\ge 4$;
since the center of $\Spin{4k}$ is ${\bf Z}_2\times {\bf Z}_2$,
we have $b_1(G(v_1,v_2);F)\le 2$. In either case, 
$b_1(G_{v_1}v_2;F)>b_1(G(v_1,v_2);F)$ which implies by
Proposition~\ref{prop:taut-sum} that $\rho_1\oplus\rho_2$ is not taut. 
\end{exs}

Recall that the \emph{slice representation} of a 
representation $\rho:G\to \mathbf O(V)$ at a point $p\in V$ 
is the representation induced by the isotropy $G_p$ 
on the normal space to the orbit $Gp$ at $p$.
The following result often works as a kind of induction.

\begin{prop}[\cite{GTh3}]\label{prop:slice}
Let $\rho:G\to \mathbf O(V)$ be a taut representation
of a compact connected Lie group $G$.
Then the slice representation of $\rho$ 
at any $p\in V$ is taut. 
\end{prop}

We now discuss a reduction principle
which in many cases considerably simplifies 
the problem of deciding whether a representation 
is taut or not. 
Let $\rho:G\to \mathbf O(V)$ be a representation of a compact
Lie group $G$ which is not assumed to be connected. Denote by $H$ a fixed
principal isotropy subgroup of the $G$-action on $V$ and let $V^H$ be the
subspace of $V$ that is left pointwise fixed by the action of
$H$. Let $N$ be the normalizer of $H$ in $G$.
Then the group $\bar N=N/H$ acts on $V^H$ with trivial principal
isotropy subgroup. Moreover, the following result is known
(\cite{G-S,Luna,LR,Schwartz,SS,Straume1}):

\begin{thm}[Luna-Richardson]\label{thm:lr}
The inclusion $V^H\to V$ induces a stratification preserving
homeomorphism between orbit spaces 
\begin{equation*}
V^H/\bar N\to V/G.
\end{equation*}
\end{thm}

The relation to tautness is expressed by the 
following result.

\begin{prop}[\cite{GTh3}]\label{prop:suf-taut}
Suppose there is a subgroup $L\subset H$ which is a 
finitely iterated $\mathbf Z_2$-extension of the identity and such that 
the fixed point sets $V^L=V^H$. Suppose also
that the reduced representation 
$\bar\rho:\bar N^0\to \mathbf O(V^H)$ is $\mathbf Z_2$-taut,
where $\bar N^0$ denotes the connected component of the identity of
$\bar N$. Then $\rho:G\to \mathbf O(V)$ is $\mathbf Z_2$-taut.
\end{prop}

We close this section with some very useful remarks. 

\begin{rmk}\label{rmk:useful}
\em
\begin{enumerate}
\item It follows from the discussion of Kuiper 
in~\cite{Ku} that \emph{if $M$ is a taut substantial submanifold
of an Euclidean space, then there exists $p\in M$ such that
the image of the second fundamental form of $M$ at $p$ 
spans the normal space of $M$ at $p$}. As a corollary, 
the codimension of $M$ is at most $n(n+1)/2$, where 
$n=\dim M$. 
\item One
defines a submanifold of an Euclidean space
to be \emph{F-tight}, or simply \emph{tight},
similarly as was done for tautness,
except that one replaces distance functions 
by height functions $h_\xi(x)=\langle x,\xi\rangle$, 
$\xi$ a nonzero vector. It turns out that tightness is invariant
under linear transformations, and a taut sumanifold
of an Euclidean space is tight. Moreover, a tight submanifold 
of an Euclidean space which is contained in a round sphere 
is taut, and in this situation the set of critical points
of a distance function will also occur as the set of 
critical points of a height function (see~\cite{CR,PT2}). 
\item Ozawa proved in~\cite{Oz} that the set of critical points 
of a distance function of a taut submanifold decomposes
into critical submanifolds which are nondegenerate 
in the sense of Bott; it follows that the so called
Morse-Bott inequalities are equalities for the function
restricted to any sublevel set; namely, the number of critical points
of the function is equal to the sum of the 
Betti numbers of the critical submanifolds, see~\cite{Bo2}.  
  
\end{enumerate}
\end{rmk}

\section{The classification}

Let $\rho:G\to\mathbf{O}(V)$ be a taut reducible
representation where $G$ is a compact connected
simple Lie group. 
Of course we may assume that
$\rho$ does not contain trivial summands. 
Write $\rho=\rho_1\oplus\rho_2$, where 
$\rho_1$ is irreducible. Then $\rho_1$ is 
orbit equivalent to the isotropy representation
of an irreducible symmetric space.
We first prove a lemma for later use,
and then we shall run through all
the possibilities for $G$ and $\rho_1$, where we find it convenient 
to consider separately the cases $G=\Spin n$ and
$G=\SO n$.

\begin{lem}\label{lem:tori}
The following representations are not taut:
\begin{enumerate}
\item $S^1\times S^1\to\U1\times\U1\times\U1$ given 
by $(e^{i\alpha},e^{i\beta})\mapsto(e^{i\alpha},e^{i\beta},
e^{i(\alpha+\beta)})$.  
\item $\SP1\times\SP1\to\SO4\times\SO4\times\SO4$ given 
by $(p,q)\mapsto(l_p,r_{\bar q},l_pr_{\bar q})$, where
$l_x$ (resp.~$r_x$) denotes left (resp.~right) translation by 
the unit quaternion $x$.
\end{enumerate}
\end{lem}

\Pf We will prove (a); assertion (b) is similar. 
Let $M$ denote the orbit through 
$p=(1,1,1)\in\CC\oplus\CC\oplus\CC$. We will
show that $M$ is not taut by exhibiting a height function
which is not perfect, see Remark~\ref{rmk:useful}(b).
The normal space $\nu_p M$ is easily seen to be spanned over $\R$ 
by $(1,0,0)$, $(0,1,0)$, $(0,0,1)$ and $(i,i,-i)$. 
Let $h:M\to\R$ be the height function defined by $p$. 
Note that $gp$, $g\in S^1\times S^1$, is a critical point
of $h$ if and only if $p\in\nu_{gp}M$, or, equivalently,
$g^{-1}p\in\nu_pM$. One immediately computes that 
$g=(\pm1,\pm1)$ or $(-\frac12\pm i\frac{\sqrt3}2,-\frac12\pm i\frac{\sqrt3}2)$,
so there are $6$ critical points. Since $M$ is a $2$-torus,
$h$ is not perfect. \EPf

\subsection{The case $G=\SO n$, $n=3$ or $n\geq5$}\label{sec:son}

Here $\rho_1$ is one 
of the following:
\begin{enumerate}
\item the vector representation on $\R^n$;
\item the adjoint representation on $\Lambda^2\R^n$, where $n\geq5$;
\item the representation on the space of traceless symmetric 
matrices $S^2_0\R^n$;
\end{enumerate}
The possibilities (b) and (c) are ruled out by Examples~\ref{exs:exs}.
Now possibility (a) is taken care of by the following 
proposition (compare~\cite{TTh1}, Examples~3.14).

\begin{prop}\label{prop:vector-sum-son}
Assume that $n\geq3$ and $\rho$ is the sum of $k>1$ 
copies of the vector representation. Then 
$\rho$ is taut. 
\end{prop}

\Pf Let $V=\R^n\oplus\cdots\oplus\R^n$, $k$ copies.
Suppose first that $k\leq n$. Let $\{e_1,\ldots,e_n\}$ be 
the canonical basis of $\R^n$, and let $p=(e_1,\ldots,e_k)\in V$. 
View $V$ as the space of real $n\times k$-matrices, and let 
$\hat G=\SO n\times\SO k$ act on $V$ by $(A,B)\cdot X=AXB^{-1}$,
where  $(A,B)\in\hat G$ and $X\in V$. Then $\hat Gp=Gp$. 
Since $(\hat G,V)$ is the isotropy representation
of the Grassmann manifold $G_k(\R^{n+k})$, we have that $Gp$
is taut. Next suppose that $k>n$ and 
let $q=(v_1,\ldots,v_k)\in V$ be an arbitrary
nonzero point. Then there is a nonsingular $k\times k$ matrix $M$ such that
right-multiplying $q$ by $M$ gives $qM=(e_1,\ldots,e_l,0,\ldots,0)\in V$,
where $1\leq l\leq n$. 
It follows from the above that $G(qM)=(Gq)M$ is taut. 
Since a taut submanifold in Euclidean space is tight,
and tightness is invariant under linear transformations, 
$Gq$ is tight. But $Gq$ lies in a sphere, and so it is taut. 
This completes the proof that $\rho$ is taut. \EPf
 
\subsection{The case $G=\SU n$, $n=3$ or $n\geq5$}\label{sec:sun}

Here $\rho_1$ is one of the following:
\begin{enumerate}
\item the vector representation on $\CC^n$;
\item the adjoint representation on $\mathfrak{su}(n)$;
\item a real form of the representation of $\SU8$ on $\Lambda^4\CC^8$.
\item the representation on the space of skew-symmetric 
matrices $\Lambda^2\CC^{2p+1}$, where $p\geq2$;
\end{enumerate}
The possibilities (b) (even if $n=4$) 
and (c) are ruled out by Examples~\ref{exs:exs}.
Consider the possibility (d). Here a principal isotropy 
subgroup $H$ is given by $p$ diagonal blocks, each isomorphic 
to $\SU2$. Denote the representation spaces of $\rho_1$ and 
$\rho_2$ by $V_1$ and $V_2$. Now there exists $v_1\in V_1$ 
such that $G_{v_1}=H\cong{\SU2}^p$. We can assume that $\rho_2$
is irreducible. If $\rho_2$ 
is the vector representation, then we can find $v_2\in V_2$
such that $Hv_2\approx S^3\times\cdots\times S^3$, $p$ factors. 
In this case $G(v_1,v_2)\approx\SU{2p+1}$; since the 
third Betti number of a compact connected simple 
Lie group is $1$, $\rho$ cannot be taut by Proposition~\ref{prop:taut-sum}. 
If $\rho_2$ is also as in (d), it is not difficult to see
that $\rho_2|G_{v_1}$
contains as a summand a representation equivalent to 
that in Lemma~\ref{lem:tori}(b), and thus 
$\rho$ cannot be taut by Proposition~\ref{prop:taut-sum}.
This shows that $\rho$ is not taut if $\rho_1$ is as 
in (d). Now (a) is covered by the following proposition. 

\begin{prop}\label{prop:vector-sum-sun}
Assume that $n\geq3$ and $\rho$ is the sum of $k$ 
copies of the vector representation. If $1\leq k<n$, then 
$\rho$ is taut. If $k\geq n$, then $\rho$ is not taut.
\end{prop}

\Pf In the case $1\leq k<n$, we need to know that the isotropy 
representation of the Grassmann manifold $G_k(\CC^{n+k})$ is 
$\mathbf{S}(\U n\times\U k)$
acting on the space of complex $n\times k$ matrices,
and it is orbit equivalent to its restriction to 
the subgroup $\SU n\times\SU k$ if $k\neq n$ (see~\cite{EH1}).
It follows as in Proposition~\ref{prop:vector-sum-son} that $\rho$ is taut. 
In the case $k\geq n$, it is enough to consider $k=n$.
Let $\{e_1,\ldots,e_n\}$ be 
the canonical basis of $\CC^n$.
The isotropy subgroup at $e_1$
is isomorphic to $\SU{n-1}$, and the slice representation 
at $e_1$ decomposes into a sum of trivial representations
and $\CC^{n-1}\oplus\cdots\oplus\CC^{n-1}$, 
$n-1$ copies. 
We use Proposition~\ref{prop:slice}
and induction to reduce to the case of 
$\SU3$ acting on $\CC^3\oplus\CC^3\oplus\CC^3$.  
Let $p=(e_1,e_2,e_3)$, and denote by $M$ the $\SU3$-orbit through
$p$. Then $M$ is the standard inclusion of $\SU3$ into the
space $M(3,\CC)$ of complex $3\times3$-matrices. 
The tangent space $T_pM$ is the Lie algebra $\mathfrak{su}(3)$,
and the normal space $\nu_pM$ is $\CC p\oplus i\mathfrak{su}(3)$.
By Remark~\ref{rmk:useful}(b), it suffices to show that
a height function is not perfect. 
Let $h:M\to\R$ be the height function defined 
by $p$. We find the critical points of $h$.  
Note that $gp$, for $g\in\SU3$, is a critical point of $h$ if and only if 
$p\in\nu_{gp}M$, or, what amounts to the same, $g^{-1}p\in\nu_pM$. Now
it is easy to see that $gp$ is a critical point of $h$
if and only if $g=\omega I$, where $\omega$ is a cubic root 
of unity and $I$ is the identity matrix, or $g$ is conjugate
to a diagonal matrix
with entries $-1$, $-1$ and $1$. It follows that
the critical set of $h$ consists of $3$ isolated
points and a submanifold diffeomorphic to $\CC P^2$, 
whence the sum of its Betti numbers  
is $6$. Since $\SU3$ has the 
homology of $S^3\times S^5$, $h$ is not perfect in the 
sense of Bott, see Remark~\ref{rmk:useful}(c). \EPf

\subsection{The case $G=\SP n$, $n\geq3$}

Here $\rho_1$ is one of the following:
\begin{enumerate}
\item the vector representation on $\CC^{2n}$;
\item the adjoint representation on $\mathfrak{sp}(n)$;
\item a real form of the $42$-dimensional representation $\Civ{}{}{}1$
of $\SP4$.
\item a real form of the representation $\Lambda^2\CC^{2n}-\CC$;
\end{enumerate}
The possibilities (b) (even if $n=2$) 
and (c) are ruled out by Examples~\ref{exs:exs}.
Consider the possibility (d). Here a principal isotropy 
subgroup $H$ is given by the diagonal embedding of $\SP1^n$ into $\SP n$,
so there exists $v_1\in V_1$ such that $G_{v_1}=H\cong{\SP1}^n$.
We can assume that $\rho_2$ is irreducible, and then $\rho_2$ is as 
in (a) or in (d). If $\rho_2$ is as in (a), 
the proof follows as in section~\ref{sec:sun} to deduce that $\rho$ is not 
taut. If $\rho_2$ is as in (d), Proposition~\ref{prop:dd} below 
implies that $\rho$ is not taut. 

\begin{prop}\label{prop:dd}
Let $V_n$ denote a real form of the representation
$\Lambda^2\CC^{2n}-\CC$ of $\SP n$, where 
$n\geq3$. Then $(\SP n, V_n\oplus V_n)$ is not taut.
\end{prop}

We postpone the proof of Proposition~\ref{prop:dd} to the end of the paper
since the methods used to prove it better belong there. 
Finally, (a) is covered by

\begin{prop}\label{prop:vector-sum-spn}
Assume that $n\geq1$ and $\rho$ is the sum of $k>1$ 
copies of the vector representation. Then 
$\rho$ is taut. 
\end{prop}

\Pf The proof is analogous to the proof of 
Proposition~\ref{prop:vector-sum-son}. \EPf

\subsection{The case $G$ is exceptional}\label{sec:exceptional}

First note that no summand of $\rho$ can be the adjoint 
representation by Example~\ref{exs:exs}(ii). 

If $G=\G$, then $\rho$ is the sum of $k$ copies of 
the $7$-dimensional representation. If $k=2$, $\rho$ is orbit
equivalent to $(\SO7,\R^7\oplus\R^7)$ (which is taut). 
If $k=3$, $\rho$ is not taut because a principal
orbit is diffeomorphic to $\G$ and thus has the 
homology of $S^3\times S^{11}$, but 
an application of Proposition~\ref{prop:taut-sum} 
would require it to have the 
homology of $S^6\times S^5\times S^3$ in case it was taut. 

If $G=\F$, then $\rho$ is the sum of $k$ copies of 
the $26$-dimensional representation. Suppose $k=2$,
$\rho=\rho_1\oplus\rho_2$. 
Then there is an isotropy subgroup $H$ of $\rho_1$ 
isomorphic to $\Spin9$. Now $\rho_2|H$ decomposes
as $\R\oplus\R^9\oplus\R^{16}$, and it is not taut
by Proposition~\ref{prop:916}. Hence $\rho$ is not taut
by Proposition~\ref{prop:taut-sum}.

$\E6$, $\E7$ and $\E8$ do not admit representations
orbit equivalent to the isotropy representation of a 
symmetric space. 

\subsection{The case $G=\Spin n$, $n=3$ or $n\geq5$}

This is case is more involved than the previous ones. 
In view of section~\ref{sec:son}, we may assume that
a summand of $\rho$ is a spin representation. Now
the only values of $n$ which need to be considered are 
$3$, $5$, $6$, $7$, $8$, $9$, $10$
and $16$. 

\subsubsection{$G=\Spin3$} Here $G=\SU2=\SP1$.
By the discussion in section~\ref{sec:son}, the 
admissible summands of $\rho$ are the vector representation
of $\SU2$ on $\CC^2$ and the representation on $\R^3$ given by $\SU2\to\SO3$.
The sum of an arbitrary number of copies of $\CC^2$ is taut by 
Proposition~\ref{prop:vector-sum-spn}. On the other hand,
$\CC^2\oplus\R^3$ is not taut, because the principal
orbit through a point $(a,b)\in\CC^2\oplus\R^3$
with $a$, $b\neq0$ is substantial and diffeomorphic to $S^3$, 
but, as mentioned in the introduction, a sphere can be taut
only in substantial codimension one. 

\subsubsection{$G=\Spin5$} Here $G=\SP2$. 
By the discussion in section~\ref{sec:son}, the 
admissible summands of $\rho$ are the vector representation
of $\SP2$ on $\CC^4$ and the representation on $\R^5$ given by 
$\SP2\to\SO5$. The situation in which $\R^5$ is not 
present is covered by Proposition~\ref{prop:vector-sum-spn}. 
On the other hand, we have

\begin{prop}\label{prop:sp2}
$\CC^4\oplus\R^5$ is not taut.
\end{prop}

\Pf Note that the principal orbits are substantial 
embeddings of $\SP2$ in $S^{12}$.
We will show that $\SP2$ can admit a taut substantial embedding 
of codimension $2$ in a sphere $S^N$ only if $N=15$ following 
an argument which appeared in~\cite{Gal}, p.~75. 

So suppose that $X$ is diffeomorphic to $\SP2$ and tautly
embedded in $S^N$ with $N\geq12$. Let $Y$ be a sufficiently 
small tubular neighborhood of $X$ in $S^N$. 
$X$ has the homology of $S^3\times S^7$, so its homology
groups vanish except in dimensions $0$, $3$, $7$ and 
$10$. Since $2\times 3\neq7$, it follows as in Proposition~2.2
of~\cite{Oz} that $Y$ is a compact proper Dupin hypersurface.
Moreover, a Morse distance function on $Y$ can have critical points
of index $0$, $3$, $7$ and $10$ only. By the Morse 
index theorem, the multiplicities of the first three
principal curvatures of $Y$ are $m_1=3$, $m_2=4$ and
$m_3=3$. According to Theorem~C in~\cite{GH}, there
exists at most $2$ different multiplicities $k$, $l$,
and $g=2$ or $4$ in case $k\neq l$. Therefore the fourth
principal curvature of $Y$ has multiplicity $m_4=4$. 
It follows that $\dim Y=14$, and hence, $N=15$. \EPf

\subsubsection{$G=\Spin6$} Here $G=\SU4$. By the discussion
in section~\ref{sec:son}, the admissible summands of $\rho$ are 
the vector representation of $\SU4$ on $\CC^4$ and the 
representation on $\R^6$ given by $\SU4\to\SO6$. 
The situation in which $\R^6$ is not 
present is covered by Proposition~\ref{prop:vector-sum-sun}. 
Also, $\CC^4\oplus\R^6$ is taut because the singular orbits 
are round spheres in $\CC^4$ and $\R^6$, and the principal
orbits are products of those. The following two 
propositions settle down this case. 

\begin{prop}
$\CC^4\oplus\R^6\oplus\R^6$ is not taut.
\end{prop}

\Pf Let $p\in\R^6$. Then the slice representation at $p$
is $\Spin5=\SP2$ acting on $\R\oplus\R\oplus\CC^4\oplus\R^5$. 
The result follows from Propositions~\ref{prop:sp2} 
and~\ref{prop:slice}. \EPf

\begin{prop}\label{prop:c4c4r6}
$\CC^4\oplus\CC^4\oplus\R^6$ is not taut.
\end{prop}

\Pf We will show that a certain orbit is not taut
by finding an explicit height function which is not perfect.
We need to have a good parametrization of the orbits.  
It is useful to use Cayley numbers. 
Recall that the Cayley algebra can be viewed
as $\Ca=\Q\oplus\Q e$ via the Cayley-Dickson process, where
$\Q=\R\langle1,i,j,k\rangle$ is the quaternion algebra
(see appendix IV.A in~\cite{Ha-La}). 
Then $\Ca=\R\langle1,i,j,k,e,ie,je,ke\rangle$. 
According to~\cite{Ch-Ri}, upon identifying $\Ca\cong\R^8$ and using
Cayley multiplication, 
\[ \Spin8=\{(A,B,C)\in\SO8\times\SO8\times\SO8:A(\xi\eta)=B(\xi)C(\eta),
\;\mbox{for all $\xi$, $\eta\in\Ca$}\}, \]
\begin{eqnarray*}
 \Spin7 & = & \{(A,B,C)\in\Spin8:A(1)=1\}, \\
        & = & \{(A,B,C)\in\Spin8:C=\tilde B\}, 
\end{eqnarray*}
where $\tilde B(x)=\overline{B(\bar x)}$, and
\[ \Spin6  =  \{(A,B,\tilde B)\in\Spin7:A(i)=i\}. \]
Also, the isomorphism $\Spin6\to\SU4$ is given by 
$(A,B,\tilde B)\mapsto B$, and the
projection $\Spin6\to\SO6$ is given by $(A,B,\tilde B)\to A$. 
Therefore the covering $\varphi:\SU4\to\SO6$ is given by
$\varphi(g)(x)=g(x)\overline{g(1)}=g(1)\overline{g(\bar x)}$,
where $g\in\SU4$ and $x\in\R^6$.  
Here we regard $\SU4$ as the subgroup of $\SO8$ defined by 
the complex structure in $\R^8$ given by
left multiplication by the element $i$. This identifies 
$\Ca\cong\CC^4$. Now (note that $i(ke)=je$)
$\CC^4 = \CC\langle1,j,e,ke\rangle$, 
$\R^6 =\R\langle j,k,e,ie,je,ke\rangle$.

Fix the base point $p=(1,j,e)\in V=\CC^4\oplus\CC^4\oplus\R^6$. 
Let $G=\SU4$ act on $V$. Then $G_p$ is trivial. 
Let $M=Gp$, principal orbit diffeomorphic to $\SU4$. $M$ can also
be parametrized by the Stiefel manifold $St_3(\CC^4)$. 
In fact, given $(z_1,z_2,z_3)\in St_3(\CC^4)$, there is a unique 
$g\in\SU4$ such that $g^{-1}(1)=z_1$, $g^{-1}(j)=z_2$, and $g^{-1}(e)=z_3$. 
Then we get $g^{-1}(1,j,e)=(z_1,z_2,z_3\bar z_1)\in M$.
View $p=(1,j,e)$ as a vector in $\nu_pM$, and let $h:M\to\R$ be the
height function defined by $p$. We have that $gp\in M$, $g\in\SU4$, is a 
critical point of $h$ if and only if $p\in\nu_pM$. 
It is easy to compute that the normal space to $M$ at $p=(1,j,e)$
is spanned by
\[ (1,0,0),\; (0,j,0),\; (0,0,e),\; 
     (j,1,0),\; (k,-i,0),\;  (je,e,j),\;
     (ke,-ie,k). \]
Now the condition that $g^{-1}p\in\nu_{p}M$ is that there exist
$A$, $B$, $C$, $D$, $E$, $F$, $G\in\R$ such that 
\[ (z_1,z_2,z_3\bar z_1) = (A+Dj+Ek+Fje+Gke,D-Ei+Bj+Fe-Gie,Fj+Gk+Ce). \]
The relations $(z_i,z_j)=\delta_{ij}$,
where $(\cdot,\cdot)$ denotes the Hermitian inner product in $\CC^4$,  
yield the following relations:
\begin{eqnarray*}
 (A+B)(D+Ei)&=&0, \\
 (F-Gi)(AB+BC+AC-F^2-G^2-D^2-E^2) & = & 0, \\
 A^2 + D^2 + E^2 + F^2 + G^2 & = & 1, \\
 A^2 - B^2 & = & 0, \\
 C^2 + F^2 + G^2 & = & 1. 
\end{eqnarray*}
The system admits the following solutions:
\begin{itemize}
\item $A=B=-C=\pm1$, $D=E=F=G=0$; 
\item $A=B=C=\pm1$,  $D=E=F=G=0$; 
\item $A=B=C=\pm\frac12$, $D=E=0$, $F^2+G^2=\frac34$;
\item $A=-B$, $C=\pm1$, $F=G=0$, $A^2+D^2+E^2=1$.
\end{itemize}
Since $g^{-1}p\mapsto gp$ is a well defined homeomorphism of $M$,
we deduce that the critical set of $h$ consists of 
$4$ points, $2$ circles and $2$ spheres.
Now the sum of the Betti numbers of the critical manifolds
of $h$ is $12$. 
Since $\SU4$ has the homology of $S^3\times S^5\times S^7$,
$M$ is not taut. \EPf

\subsection{$G=\Spin7$}\label{sec:spin7} By the discussion
in section~\ref{sec:son}, the admissible summands of $\rho$ are 
the vector representation on $\R^7$ and the spin
representation on $\R^8$. We first note that 
$\R^8\oplus\R^7$ is taut because the singular orbits 
are round spheres in $\R^8$ and $\R^7$, and the principal
orbits are products of those. Moreover, 
$\R^8\oplus\R^8$ and $\R^8\oplus\R^8\oplus\R^8$ are taut
because $\Spin7$ is transitive on the Stiefel manifolds
$St_2(\R^8)$ and $St_3(\R^8)$, so the actions of $\Spin7$ 
on these spaces are orbit equivalent to the actions of $\SO8$. 
We also note that if $\rho$ has $4$ summands and $\R^8$ is one of them, 
say $V_1$, then the slice representation at a point in $V_1$ 
is $\G$ acting on $\R^7\oplus\R^7\oplus\R^7$, which is not taut
by the discussion in section~\ref{sec:exceptional};
hence, $\rho$ is not taut by Proposition~\ref{prop:slice}. 
We finish this case with the following two propositions. 

\begin{prop}\label{prop:r7r7r8}
$\R^7\oplus\R^7\oplus\R^8$ is taut. 
\end{prop}

\Pf We shall use the reduction principle 
as described in Proposition~\ref{prop:suf-taut}. 
In order to have a good description of the 
representation, we resort to Cayley numbers
as in the proof of
Proposition~\ref{prop:c4c4r6}.
View
$\R^8=\R\langle1,i,j,k,e,ie,je,ke\rangle$
and $\R^7=\R\langle i,j,k,e,ie,je,ke\rangle$. Let 
$G=\Spin7$, $V=\R^7\oplus\R^7\oplus\R^8$. 
The action of $G$ on $V$ is given by $(A,B,\tilde B)\mapsto(A,A,B)$. 
The isotropy of $G$ at $p=(i,j,1)\in V$ is
\begin{eqnarray*}
H & = & \{(A,A,A)\in\Spin8:\;\mbox{$A\in\SP2$, $A$ fixes $1$} \} 
\cong\SP1,
\end{eqnarray*}
where we regard $\SP2$ as the subgroup of $\SO8$ defined by 
the complex structures in $\R^8$ given by the
left multiplications by the elements $i$, $j$. 
This identifies $\R^8\cong\Q\langle1,e\rangle$. 

The description of $H$ shows that the cohomogeneity of $(G,V)$ is $4$
and the fixed point subspace
\[ V^H = \R\langle i,j,k\rangle\oplus \R\langle i,j,k\rangle\oplus
 \R\langle1,i,j,k\rangle \cong\R^{10}. \]
It follows from Theorem~\ref{thm:lr} that $\dim\bar N=6$. 
The normalizer $N$ of $H$ in $G$ is the same as 
the stabilizer of $V^H$ in $G$. Suppose that
$(A,B,\tilde B)\in N$. Then we can write
\[ A = \left(\begin{array}{cc}
             A_1 & 0 \\
             0 & A_2 \end{array}\right), \]
where $A_1$, $A_2\in\SO4$, $A_1(1)=1$, and we view 
$\R^4=\R\langle 1,i,j,k\rangle$. Since $\SP1\times\SP1\to\SO4$,
$(p,q)\mapsto l_pr_{\bar q}$ (notation as in Lemma~\ref{lem:tori})
is a double covering, we can write $A_2=l_pr_{\bar q}$ 
for unique $(p,q)$ modulo $\pm1$. Similarly,
$\SP1\to\SO3$ , $s\mapsto l_s r_{\bar s}$ is a double covering, 
so we can write $A_1=l_s r_{\bar s}$ 
for a unique $s$ modulo $\pm1$. We deduce that (compare~\cite{Ch-Ri},
section~2)
\begin{equation}\label{N}
 (A,B,\tilde B)=\left(\left(\begin{array}{cc}
                                l_sr_{\bar s} & 0 \\
                                0 & l_pr_{\bar q} 
                         \end{array}\right),
                   \left(\begin{array}{cc}
                                l_sr_{\bar q} & 0 \\
                                0 & l_pr_{\bar s} 
                         \end{array}\right),
                    \left(\begin{array}{cc}
                                l_qr_{\bar s} & 0 \\
                                0 & l_pr_{\bar s} 
                         \end{array}\right)
                  \right).
\end{equation}
Therefore $N$ consists of the elements of the form~(\ref{N})
for $p$, $q$, $s\in\SP1$, and $H$ consists of the elements
with $q=s=1$. Now  
\[\bar N=N/H\cong\SP1\times_{\mathbf{Z}_2}\SP1=\{(q,s)\in\SP1\times\SP1:(q,s)\sim(-q,-s)\}, \]
the action of $\bar N$ on $V^H$ is given by
\[ (q,s)\in\bar N\mapsto(l_sr_{\bar s},l_sr_{\bar s},l_sr_{\bar q})
\in\SO3\times\SO3\times\SO4, \]
and thus it is orbit equivalent to the product of the standard action 
of $\SO3$ on $\R^3\oplus\R^3$ by the standard action of 
$\SP1$ on $\CC^2$. Since these are taut representations,
we deduce that $(\bar N,V^H)$ is also taut. Now let
$L$ be the $\Z_2$-subgroup of $H$ generated by the 
element~(\ref{N}) with $q=s=1$, $p=-1$. Then $V^L=V^H$.
It follows from Proposition~\ref{prop:suf-taut} that
$(G,V)$ is taut. \EPf

\begin{prop}
$\R^7\oplus\R^8\oplus\R^8$ is not taut. 
\end{prop}

\Pf We use a method similar to that of the proof
of Proposition~\ref{prop:r7r7r8}. Let 
$G=\Spin7$, $V=\R^7\oplus\R^8\oplus\R^8$.
The action of $G$ on $V$ is given by
$(A,B,\tilde B)\mapsto(A,B,B)$.
The isotropy of $G$ at $p=(i,1,j)\in V$ is
\[ H=\{(A,A,A)\in\Spin8:\;\mbox{$A\in\SU4$, $A$ fixes $1$, $j$} \}\cong\SU2, \]
and the cohomogeneity of $(G,V)$ is $5$. The fixed point subspace
\[ V^H = \R\langle i,j,k\rangle\oplus \R\langle1,i,j,k\rangle\oplus
 \R\langle1,i,j,k\rangle \cong\R^{11}, \]
and $\dim\bar N=6$. Now $N$, $H$ and $\bar N$ are as in 
Proposition~\ref{prop:r7r7r8}, and the 
action of $\bar N$ on $V^H$ is given by
\[ (q,s)\in\bar N\mapsto(l_sr_{\bar s},l_sr_{\bar q},l_sr_{\bar q})
\in\SO3\times\SO4\times\SO4. \]
Let $M=Gp$, and let $h$ denote the height function
defined by $p$ on $M$. It is not difficult to see that
the critical set of the restriction $h|M\cap V^H$
coincides with the critical set of $h$ (compare Lemma~3.17 in~\cite{GTh1}). 
But $M\cap V^H=\bar Np$, and a 
tedious computation shows that the sum of the Betti 
numbers of the critical set  
of $h|\bar Np$ is $12$. If $M$ was taut, it would have to have the 
homology of $S^5\times S^6\times S^7$ by Proposition~\ref{prop:taut-sum},
so the sum of its Betti numbers would have to be $8$. 
It follows that $M$ is not taut. \EPf

\subsection{$G=\Spin8$} By the discussion
in section~\ref{sec:son}, the admissible summands of $\rho$ are 
the vector representation which we denote by $\R_0^8$, and the 
half-spin representations, which we denote by $\R^8_+$ and
$\R^8_-$. The group of automorphisms of $\Spin8$ is isomorphic to 
the dihedral group of degree $3$, and it permutes the representations
$\R^8_0$, $\R^8_+$, $\R^8_-$, so this reduces the number of cases
to be considered. We now note that $\R^8_0\oplus\R^8_+$ 
is taut because the principal orbits are products of spheres;
up to permutations, there are no other representations with two summands
which need to be considered. Similarly, in the case of three
summands, there are only two cases to be considered,
see Propositions~\ref{prop:00+} and~\ref{prop:0+-}.
In the case of four summands, at least two of them 
coincide, and we can assume that those are $\R^8_0$. 
So there are three cases to be considered:
$\R^8_0\oplus\R^8_0\oplus\R^8_+\oplus\R^8_+$,
$\R^8_0\oplus\R^8_0\oplus\R^8_+\oplus\R^8_-$, and
$\R^8_0\oplus\R^8_0\oplus\R^8_0\oplus\R^8_+$;
the first one of these is not taut since a slice
representation contains $(\Spin7,\R^7\oplus\R^8\oplus\R^8)$,
which is not taut, and we can apply Proposition~\ref{prop:slice};
the second one is not taut because it contains 
$(\Spin8,\R^8_0\oplus\R^8_+\oplus\R^8_-)$, which is not taut
by Proposition~\ref{prop:0+-}; and the third one is taut
by Proposition~\ref{prop:000+}. In the case 
of five summands, there is always a slice representation
equivalent to $(\G,\R^7\oplus\R^7\oplus\R^7)$, which is 
not taut, and we can apply Proposition~\ref{prop:slice}. 

\begin{prop}\label{prop:00+}
$\R^8_0\oplus\R^8_0\oplus\R^8_+$ is taut. 
\end{prop}

\Pf We use a method similar to that of the proof
of Proposition~\ref{prop:r7r7r8}. Let 
$G=\Spin8$, $V=\R^8_0\oplus\R^8_0\oplus\R^8_+$.
The action of $G$ on $V$ is given by
$(A,B,C)\mapsto(A,A,B)$. The isotropy of $G$ at $p=(1,i,1)\in V$ is
\[ H=\{(A,A,A)\in\Spin8:\;\mbox{$A\in\SU4$ fixes $1$} \}\cong\SU3, \]
and the cohomogeneity of $(G,V)$ is $4$.
The fixed point subspace
\[ V^H = \R\langle 1,i\rangle\oplus \R\langle1,i\rangle\oplus
 \R\langle1,i\rangle \cong\R^6, \]
and $\dim\bar N=2$. We now construct two one-parameter
subgroups of $N$ which do not lie in $H$. Let $A\in\SO8$ be the rotation
by $\theta$ on the plane $\R\langle1,i\rangle$ fixing its
orthogonal complement, and let $B(x)=e^{\frac{i\theta}2}x$,
$C(x)=xe^{\frac{i\theta}2}$, for $x\in\Ca$. Then 
$(A,B,C)\in N$. We denote this transformation by $t_\theta$. 
Next, let $A\in\SO8$ fix $1$, $i$, and let $B\in\SU4$
act on $\CC\langle 1,j,e,ke\rangle$ by the matrix 
$\mathrm{diag}(e^{i\varphi},e^{-i\varphi},1,1)$. Then 
$(A,B,\tilde B)\in N$. We denote this transformation by $s_\varphi$. 
Now 
\[ \bar N^0=N^0/H\cong S^1\times S^1=\{(t_\theta,s_\varphi)\}, \]
and the action of $\bar N^0$ on $V^H$ is given by
\[ (t_\theta,s_\varphi)\in\bar N^0\mapsto
(e^{i\theta},e^{i\theta},e^{i(\frac\theta2+\varphi)})
\in\U1\times\U1\times\U1. \]
This action is clearly taut. Let $L$ be the subgroup of $H$
generated by the diagonal matrices with $\pm1$ entries. Then 
$V^L=V^H$, and $(G,V)$ is taut by Proposition~\ref{prop:suf-taut}. \EPf

\begin{prop}\label{prop:0+-}
$\R^8_0\oplus\R^8_+\oplus\R^8_-$ is not taut. 
\end{prop}

\Pf Here the action of $G$ on $V$ is given by
$(A,B,C)\mapsto(A,B,C)$. The isotropy of $G$ at $p=(1,1,i)\in V$ is
\[ H=\{(A,A,A)\in\Spin8:\;\mbox{$A\in\SU4$ fixes $1$} \}\cong\SU3, \]
and the cohomogeneity of $(G,V)$ is $4$.
The fixed point subspace $V^H$ and $\bar N^0$ 
are as in Proposition~\ref{prop:00+}, and the 
action of $\bar N^0$ on $V^H$ is given by
\[ (t_\theta,s_\varphi)\in\bar N^0\mapsto
(e^{i\theta},e^{i(\frac\theta2+\varphi)},e^{i(\frac\theta2-\varphi)})
\in\U1\times\U1\times\U1. \]
Since this action is equivalent to that of Lemma~\ref{lem:tori}(a),
it is not taut. It follows that $(G,V)$ is not taut by the final 
argument in the proof of Lemma~6.11 in~\cite{GTh3}. \EPf

\begin{prop}\label{prop:000+}
$\R^8_0\oplus\R^8_0\oplus\R^8_0\oplus\R^8_+$ is taut.
\end{prop}

\Pf Here the action of $G$ on $V$ is given by
$(A,B,C)\mapsto(A,A,A,B)$. The isotropy of $G$ at $p=(1,i,j,1)\in V$ is
the same $H$ as in Proposition~\ref{prop:r7r7r8}, the cohomogeneity is $7$, 
the fixed point subspace
\[ V^H=\R\langle1,i,j,k\rangle\oplus \R\langle1,i,j,k\rangle\oplus
 \R\langle1,i,j,k\rangle \oplus
 \R\langle1,i,j,k\rangle \cong\R^{16}, \]
and so $\dim\bar N=9$. As in Proposition~\ref{prop:r7r7r8},
we compute that
\[ N =\left\{ \left(\left(\begin{array}{cc}
                                l_sr_{\bar t} & 0 \\
                                0 & l_pr_{\bar q} 
                         \end{array}\right),
                   \left(\begin{array}{cc}
                                l_sr_{\bar q} & 0 \\
                                0 & l_pr_{\bar t} 
                         \end{array}\right),
                    \left(\begin{array}{cc}
                                l_qr_{\bar t} & 0 \\
                                0 & l_pr_{\bar s} 
                         \end{array}\right)
                  \right):p,q,s,t\in\SP1\right\}, \]
and $(\bar N,V^H)$ is 
\[ (q,s,t)\in\bar N\mapsto(l_sr_{\bar t},l_sr_{\bar t},l_sr_{\bar t},
l_sr_{\bar q})\in\SO4\times\SO4\times\SO4\times\SO4. \]
This action is orbit equivalent to the product of
$(\SO4,\R^4\oplus\R^4\oplus\R^4)$ and $(\SP1,\CC^2)$,
hence, taut. We take $L$ as in Proposition~\ref{prop:r7r7r8} and 
we get that $(G,V)$ is taut by Proposition~\ref{prop:suf-taut}. \EPf

\subsection{G=\Spin9} By the discussion
in section~\ref{sec:son}, the admissible summands of $\rho$ are 
the vector representation on $\R^9$ and the spin
representation on $\R^{16}$. Note that
$\R^{16}\oplus\R^{16}\oplus\R^{16}$ is not taut
since a slice representation is 
$(\Spin7,\R^7\oplus\R^8\oplus\R^7\oplus\R^8)$.
The other possibilities are covered by the following two propositions.

\begin{prop}\label{prop:1616}
$\R^{16}\oplus\R^{16}$ is taut.
\end{prop}

\Pf We need to have a good description of the spin 
representation of $\Spin9$. We start by letting
 $\{e_1,\ldots,e_n\}$ be the canonical basis of $\R^n$, 
and recalling that the Clifford algebra $\mathcal{C}\ell(n)$ 
(resp.~$\mathcal{C}\ell_+(n)$)
is the real associative algebra with unit generated
by $e_1,\ldots,e_n$ subject to the relations 
$e_ie_j+e_je_i=-2\delta_{ij}$ (resp.~$e_ie_j+e_je_i=+2\delta_{ij}$).
The group $\Spin n$ (resp.~$\mathbf{Spin}_+(n)$)
is the multiplicative subgroup of $\mathcal{C}\ell(n)$
(resp.~$\mathcal{C}\ell_+(n)$) consisting of even products
of elements in the unit sphere of $\R^n$. It is clear that there is 
an isomorphism $\mathcal{C}\ell(n)\otimes\CC\to\mathcal{C}\ell_+(n)
\otimes\CC$, 
induced by $e_i\mapsto\sqrt{-1}e_i$, which restricts to an isomorphism
$\Spin n\to\mathbf{Spin}_+(n)$ (see e.g. chapters~13 and~15 in~\cite{Post}). 

Now view
\[ \R^9=\R\oplus\Ca,\quad\R^{16}=\Ca\oplus\Ca, \]
where $\Ca=\R\langle 1,e,i,j,k,ei,ej,ek \rangle$,
and write $\{e_0;e_1,\ldots,e_8\}$ for the basis
$\{1;1,\ldots,ek\}$ of $\R^9$. Define
\[ \varphi:\R^9\to M(16,\R),
\qquad(r,u)\mapsto\left(\begin{array}{cc}
                               rI_8 & R_u \\
                               R_{\bar u} & -rI_8 
                       \end{array}\right), \]
where $r\in\R$, $u\in\Ca$, and $R_u:\Ca\to\Ca$ is
right Cayley multiplication. 
Then $\varphi(r,u)^2=(r^2+||u||^2)I_{16}$.
It follows that $\varphi$ induces a 
homomorphism
$\mathcal{C}\ell_+(9)\to M(16,\R)$. Restricting to 
$\mathbf{Spin}_+(9)$ and identifying 
$\Spin9\cong\mathbf{Spin}_+(9)$, we finally get the 
spin representation $\Delta_9:\Spin9\to\SO{16}$. 

Now consider $G=\Spin9$ acting on $V=\R^{16}\oplus\R^{16}$
via $\Delta_9\oplus\Delta_9$,
where $\R^{16}=\Ca\oplus\Ca$. 
The principal isotropy subgroup $H$ at the point $((1,0),(e,1))\in V$
is isomorphic to $\SU3$, and $\Delta_9(H)$ consists of 
matrices of the form 
\begin{equation}\label{H}
 \left( \begin{array}{cccccc}
           1&&&&& \\
           &1&&&& \\
           &&\begin{array}{cc} A & B \\ -B & A \end{array}&&&\\
           &&&1&&\\
           &&&&1&\\
           &&&&&\begin{array}{cc} A & B \\ -B & A \end{array}
          \end{array} \right)\in\SO{16}, 
\end{equation}
where $A+iB\in\SU3$. Now the cohomogeneity of $(G,V)$ is $4$,
the fixed point subspace 
\[ V^H=\R\langle(1,0),(e,0),(0,1),(0,e)\rangle\oplus
\R\langle(1,0),(e,0),(0,1),(0,e)\rangle\subset\R^{16}\oplus\R^{16},\]
and 
$\dim\bar N=4$. Using the above description of $\Delta_9$, one can check
that $e_0e_1$, $e_1e_2$, $e_0e_2$ belong to $N$ and generate a
subgroup isomorphic to $\SU2$. Moreover $e_3e_4e_5e_6e_7e_8$
centralizes this subgroup and also belongs to $N$. Hence
$\bar N^0\cong\U2$, and $(\bar N^0,V^H)$ is $(\U2,\CC^2\oplus\CC^2)$;
this representation is taut by an argument similar to one used in the proof
of Proposition~\ref{prop:vector-sum-sun}, based on the fact that the isotropy 
representation of the Grassmann manifold $G_2(\CC^2)$ is orbit 
equivalent to $\U2\times\U2$ acting on complex $2\times 2$ matrices. 
Let $L$ be the
subgroup of $H$ generated by the elements~(\ref{H})
with $A$ diagonal with $\pm1$ entries and $B=0$. Then $V^L=V^H$. 
Thus, $(G,V)$ is taut by Proposition~\ref{prop:suf-taut}. \EPf

\begin{prop}\label{prop:916}
$\R^9\oplus\R^{16}$ is not taut.
\end{prop}

\Pf We use the description of the spin representation 
given in the proof of Proposition~\ref{prop:1616}.
One can check that the principal isotropy subgroup $H$ 
at $(e_0,(1,1))\in\R^9\oplus(\Ca\oplus\Ca)$ is isomorphic to $\G$,  
$V^H=\R\langle e_0,e_1\rangle\oplus(\R1\oplus\R1)\subset\R^9\oplus
(\Ca\oplus\Ca)$, the cohomogeneity is $3$, and so $\dim\bar N=1$. 
It then follows that $\theta\mapsto\cos\theta1+\sin\theta(e_0e_1)$
defines a one-parameter subgroup in $\bar N$ which acts 
on $(\R1\oplus\R1)$ as a rotation by an angle of $\theta$,
and acts on $\R\langle e_0,e_1\rangle$  
as a rotation by an angle of $2\theta$. 
Therefore $(\bar N,V^H)$ is not taut. 
It follows that $(G,V)$ is not taut by the final 
argument in the proof of Lemma~6.11 in~\cite{GTh3}. \EPf

\subsection{G=\Spin{10}} By the discussion
in section~\ref{sec:son}, the admissible summands of $\rho$ are 
the vector representation on $\R^{10}$, and the 
half-spin representations on $\CC^{16}_+$ and
$\CC^{16}_-$. It is clear that the following two propositions
cover all possibilities. 

\begin{prop}
$\R^{10}\oplus\CC^{16}_+$ is not taut.
\end{prop}

\Pf We extend the ideas of Proposition~\ref{prop:1616}.
Let $\mathcal C\ell^0(n)$ denote the ``even'' part of
$\mathcal C\ell(n)$, namely the subalgebra
of $\mathcal C\ell(n)$ consisting of even products of elements
in $\R^n$. Then $\Spin n$ is a subgroup of $\mathcal C\ell^0(n)$,
and an isomorphism $\mathcal C\ell^0(n)\cong\mathcal C\ell(n-1)$
is given by 
\[ \left\{ \begin{array}{ll}
            e_ie_j \mapsto e_ie_j, & \mbox{if $i<j<n$}, \\
            e_ie_n \mapsto e_i, & \mbox{if $i<n$}. 
           \end{array} \right. \]
View $\R^9=\R\oplus\Ca$ and $\R^{16}=\Ca\oplus\Ca$ as 
in Proposition~\ref{prop:1616}, and define 
\[ \varphi_{\pm}:\R^9\to M(16,\CC), 
\qquad(r,u)\mapsto\pm\sqrt{-1}\left(\begin{array}{cc}
                               rI_8 & R_u \\
                               R_{\bar u} & -rI_8 
                       \end{array}\right), \]
where $r\in\R$, $u\in\Ca$, and $R_u:\Ca\to\Ca$ is
right Cayley multiplication. 
Then $\varphi_{\pm}(r,u)^2=-(r^2+||u||^2)I_{16}$.
It follows that $\varphi_{\pm}$ induce  
homomorphisms $\mathcal{C}\ell(9)\to M(16,\CC)$. 
Now $\Spin{10}\subset\mathcal{C}\ell^0(10)\cong\mathcal{C}\ell(9)$,
so these homomorphisms restrict to the half-spin
representations $\Delta_{10}^{\pm}:\Spin{10}\to\U{16}$. 
Note that $\omega=e_0e_1e_2e_3e_4e_5e_6e_7e_8e_9$ belongs to
the center of $\Spin{10}$ and $\Delta_{10}^{\pm}(\omega)=\pm\sqrt{-1}I_{16}$.
It follows that $\Delta_{10}^+$ and $\Delta_{10}^-$ are not equivalent.
It is also clear that $\Delta_{10}^{\pm}|\Spin9=\Delta_9\oplus\Delta_9$. 

Next consider $G=\Spin{10}$ acting on $V=\R^{10}\oplus\CC^{16}_+$. 
We view $\CC^{16}_+=\R^{16}\oplus\sqrt{-1}\R^{16}$,
$\Spin9$-invariant decomposition, where $\R^{16}=\Ca\oplus\Ca$.
A principal isotropy subgroup can be taken to be the same subgroup $H$
as in Proposition~\ref{prop:1616}, and the fixed subspace
\begin{eqnarray*}
\lefteqn{V^H=\R\langle e_0,e_1,e_2,e_9\rangle\oplus
\R\langle(1,0),(e,0),(0,1),(0,e)\rangle}\\
&&\qquad\oplus
\R\langle(\epsilon1,0),(\epsilon e,0),(0,\epsilon1),
(0,\epsilon e)\rangle\subset\R^{10}\oplus\R^{16}\oplus\epsilon\R^{16}, 
\end{eqnarray*}
where $\epsilon=\sqrt{-1}$. Now the cohomogeneity of $(G,V)$ 
is $5$ and $\dim\bar N=7$. 

It is not difficult to see that $\bar N^0$ is locally isomorphic 
to $\U1\times\SU2_1\times\SU2_2$, where the $\U1$-factor is generated
by $e_3e_4e_5e_6e_7e_8$ 
and the Lie algebras of the $\SU2$-factors are respectively
spanned by $e_0e_1+e_2e_9$, $e_0e_2-e_1e_9$, $e_0e_9+e_1e_2$ and
$e_0e_1-e_2e_9$, $e_0e_2+e_1e_9$, $e_0e_9-e_1e_2$. We want to describe the
action of $\bar N^0$ on $V^H$. For that purpose, it is convenient 
to set $\R^4=V^H\cap\R^{10}$ and 
$\CC^4=V^H\cap\CC^{16}_+$. Then it can be shown 
that there is a decomposition $\CC^4=\CC^2_1\oplus\CC^2_2$
such that $\SU2_1\times\SU2_2$ acts by the product
of the standard representations on $\CC^2_1\oplus\CC^2_2$ and
it acts on $\R^4$ by $\SU2_1\times\SU2_2\to\SO4$. 
Moreover, $\U1$ acts scalarly on $\CC^2_1$, $\CC^2_2$,
and trivially on $\R^4$. 
We finally get that $(\bar N^0,V^H)$ is equivalent to 
\[ (e^{j\theta},p,q)\in(\U1\times\SP1\times\SP1)/\Z_2
\mapsto(l_pr_{e^{-j\theta}},l_qr_{e^{j\theta}},l_pr_{\bar q})\in
\SO4\times\SO4\times\SO4, \]
where we have identified $V^H=\Q\oplus\Q\oplus\Q$. 
It is also important to note that $\bar N$ is not connected,
and the element $e_1e_5e_7e_6$ lies in $\bar N\setminus\bar N^0$. 

Finally, consider the $\bar N$-orbit of 
$x=(1,1,1)\in\Q\oplus\Q\oplus\Q$, and let $h$ be the height function defined
by $x$. A careful calculation shows that the sum of the Betti
numbers of the critical set of $h$ on $\bar N^0x$
is $12$. Therefore, on $\bar Nx$, this sum is at least $24$. 
The critical set of $h$ 
on $\bar Nx$ is the same as its critical set on $M=Gx$. 
If $M$ is taut, it has the homology of 
$S^{15}\times S^9\times S^7\times S^6$ by 
Proposition~\ref{prop:taut-sum}, so the sum of the Betti 
numbers of $M$ has to be $16$. Hence, $M$ is not taut. \EPf

\begin{prop}
$\CC^{16}_+\oplus\CC^{16}_+$ and $\CC^{16}_+\oplus\CC^{16}_-$ 
are not taut.
\end{prop}

\Pf In both representations, the principal isotropy subgroup of the first
summand acts on the second summand by a representation that contains
a summand equivalent to $(\SU4,\CC^4\oplus\R^6\oplus\R^6\oplus\CC^4)$,
which is not taut. Hence we can apply Proposition~\ref{prop:taut-sum}. \EPf

\subsection{G=\Spin{16}} This case is ruled out because 
the spin representation on 
$\R^{128}$ cannot be a summand of a taut representation of
$\Spin{16}$ by the argument of Example~\ref{exs:exs}(i). 

\bigskip

\textit{Proof of Proposition~\ref{prop:dd}.}
Consider first the representation $(\SP n,V_n)$.
Let $K$ be $\SP1\times\SP{n-1}$ diagonally
embedded into $\SP n$. Then there exists a point
in $V_n$ whose isotropy subgroup is $K$, and such that 
its slice representation contains as a summand $V_{n-1}$. 
This implies that $V_n\oplus V_n$ admits a slice
representation containing $V_{n-1}\oplus V_{n-1}$. 
By Proposition~\ref{prop:slice} and induction on $n$,
it is now enough to prove that $(\SP3,V_3\oplus V_3)$
is not taut. 

The principal isotropy subgroup of $(\SP3,V_3)$ is the diagonal
embedding of $\SP1^3$ into $\SP3$; call it $K_1$. 
Now $V_3$, considered as a representation of $K_1$,
decomposes into two copies of the trivial representation
and a representation $W$ which, upon identification with 
$\Q\oplus\Q\oplus\Q$, is orbit equivalent to (notation as in 
Lemma~\ref{lem:tori})
\[ (p,q,s)\in\SP1^3\mapsto(l_pr_{\bar q}, l_pr_{\bar s}, l_qr_{\bar s})
\in\SO4\times\SO4\times\SO4. \]
By Proposition~\ref{prop:taut-sum}, it is enough to show that
$(K_1,W)$ is not taut, and, for that purpose, we will
apply the reduction principle described in Proposition~\ref{prop:suf-taut}
to $(K_1,W)$.  

The principal isotropy subgroup of $(K_1,W)$
at the point $(1,i,j)\in\Q\oplus\Q\oplus\Q$ is the circle
subgroup $H=\{(e^{kt},e^{kt},e^{-kt}):t\in\R\}$ of $K_1$. 
Therefore the cohomogeneity of $(K_1,W)$ 
is $4$, the fixed point subspace of $H$
is $W^H=\R\langle 1,k\rangle\oplus\R\langle i,j\rangle\oplus
\R\langle i,j\rangle$, and so the dimension of the
normalizer $N$ of $H$ in $K_1$ is $3$. It is clear 
that $N^0=\{(e^{ka},e^{kb},e^{kc}):a,b,c\in\R\}$. 
Consider the one-parameter subgroups of $N$ given
by $\varphi_a=(e^{ka},e^{-ka},1)$ and $\psi_b=(1,e^{kb},e^{kb})$.
Then $\varphi_a$ and $\psi_b$ generate $\bar N^0$,
and $(\bar N^0,V^H)$ is 
$(\varphi_a,\psi_b)\mapsto(e^{k(2a-b)},e^{k(a+b)},e^{k(-a+2b)})$,
which is not taut by Lemma~\ref{lem:tori}(a). 
It follows that $(K_1,W)$ is not taut by the final 
argument in the proof of Lemma~6.11 in~\cite{GTh3}. \EPf

\providecommand{\bysame}{\leavevmode\hbox to3em{\hrulefill}\thinspace}


\begin{thebibliography}{HPT88}

\bibitem[Ban70]{Ba2}
T.~F. Banchoff, \emph{The spherical two-piece property and tight surfaces in
  spheres}, {J. Differential Geom.} \textbf{4} (1970), 193--205.

\bibitem[Bot54]{Bo2}
R.~Bott, \emph{Nondegenerate critical manifolds}, Ann. of Math. (2) \textbf{60}
  (1954), 248--261.

\bibitem[BS58]{BS}
R.~Bott and H.~Samelson, \emph{Applications of the theory of {M}orse to
  symmetric spaces}, {Amer. J. Math.} \textbf{80} (1958), 964--1029, Correction
  in {Amer. J. Math. {\bf83} (1961), 207--208}.

\bibitem[CL57]{CL1}
S.~S. Chern and R.~Lashof, \emph{On the total curvature of immersed manifolds},
  {Amer. J. Math.} \textbf{79} (1957), 306--318.

\bibitem[CR78]{CR2}
T.~E. Cecil and P.~J. Ryan, \emph{Focal sets, taut embeddings and the cyclides
  of {D}upin}, {Math. Ann.} \textbf{236} (1978), 177--190.

\bibitem[CR85]{CR}
T.~E. Cecil and P.~J. Ryan, \emph{Tight and taut immersions of manifolds},
  Research Notes in Mathematics, no. 107, Pitman, 1985.

\bibitem[CR98]{Ch-Ri}
L.~M. Chaves and A.~Rigas, \emph{From the triality viewpoint}, {Note Mat.}
  \textbf{18} (1998), no.~2, 155--163.

\bibitem[CW72]{CW1}
S.~Carter and A.~West, \emph{Tight and taut immersions}, {Proc. London. Math.
  Soc.} \textbf{25} (1972), 701--720.

\bibitem[EH99]{EH1}
J.~Eschenburg and E.~Heintze, \emph{On the classification of polar
  representations}, {Math. Z.} \textbf{232} (1999), 391--398.

\bibitem[Gal93]{Gal}
B.~Galemann, \emph{Tautness and linear representations of the classical compact
  groups}, Ph.D. thesis, University of Notre Dame, 1993.

\bibitem[GH91]{GH}
K.~Grove and S.~Halperin, \emph{Elliptic isometries, condition {(C)} and proper
  maps}, {Arch. Math. (Basel)} \textbf{56} (1991), 288--299.

\bibitem[GS00]{G-S}
K.~Grove and C.~Searle, \emph{Global {$G$}-manifold reductions and
  resolutions}, {Ann. Global Anal. and Geom.} \textbf{18} (2000), 437--446,
  Special issue in memory of Alfred Gray (1939--1998).

\bibitem[GT]{GTh1}
C.~Gorodski and G.~Thorbergsson, \emph{Representations of compact {L}ie groups
  and the osculating spaces of their orbits}, preprint, Univ. of Cologne, 2000
  (also E-print math.~DG/0203196).

\bibitem[GT02]{GTh2}
C.~Gorodski and G.~Thorbergsson, \emph{Cycles of {B}ott-{S}amelson type for
  taut representations}, {Ann. Global Anal. Geom.} \textbf{21} (2002),
  287--302.

\bibitem[GT03]{GTh3}
C.~Gorodski and G.~Thorbergsson, \emph{The classification of taut irreducible
  representations}, {J. Reine Angew. Math.} \textbf{555} (2003), 187--235.

\bibitem[Heb88]{He}
J.~Hebda, \emph{The possible cohomology ring of certain types of taut
  submanifolds}, {Nagoya Math. J.} \textbf{111} (1988), 85--97.

\bibitem[HL82]{Ha-La}
R.~Harvey and H.~B. {Lawson Jr.}, \emph{Calibrated geometries}, {Acta Math.}
  \textbf{148} (1982), 47--157.

\bibitem[HPT88]{HPT}
W.-Y. Hsiang, R.~S. Palais, and C.-L. Terng, \emph{The topology of
  isoparametric submanifolds}, {J. Differential Geom.} \textbf{27} (1988),
  423--460.

\bibitem[Kui61]{Ku}
N.~H. Kuiper, \emph{Sur les immersions \`a courbure totale minimale}, S\'
  eminaire de Topologie et G\' eometrie Diff\'erentielle C. Ereshmann, Paris,
  vol.~II, 1961, Recueil d'expos\'es faits en 1958-1959-1960.

\bibitem[LR79]{LR}
D.~Luna and R.~W. Richardson, \emph{A generalization of the {C}hevalley
  restriction theorem}, {Duke Math. J.} \textbf{46} (1979), 487--496.

\bibitem[Lun75]{Luna}
D.~Luna, \emph{Adh\'erences d'orbite et invariants}, {Invent. Math.}
  \textbf{29} (1975), 231--238.

\bibitem[Oza86]{Oz}
T.~Ozawa, \emph{On the critical sets of distance functions to a taut
  submanifold}, {Math.~Ann.} \textbf{276} (1986), 91--96.

\bibitem[Pos86]{Post}
M.~Postnikov, \emph{{L}ie groups and {L}ie algebras, {L}ectures in geometry,
  {S}emester {V}}, Mir, Moscow, 1986, Translated from the Russian by Vladimir
  Shokurov.

\bibitem[PT88]{PT2}
R.~S. Palais and C.-L. Terng, \emph{Critical point theory and submanifold
  geometry}, Lect. Notes in Math., no. 1353, Springer-Verlag, 1988.

\bibitem[PT89]{PTh2}
U.~Pinkall and G.~Thorbergsson, \emph{Taut $3$-manifolds}, {Topology}
  \textbf{28} (1989), 389--401.

\bibitem[Sch80]{Schwartz}
G.~W. Schwartz, \emph{Lifting smooth homotopies of orbit spaces}, {I.H.E.S.
  Publ. in Math.} \textbf{51} (1980), 37--135.

\bibitem[SS95]{SS}
T.~Skjelbred and E.~Straume, \emph{A note on the reduction principle for
  compact transformation groups}, preprint, 1995.

\bibitem[Str94]{Straume1}
E.~Straume, \emph{On the invariant theory and geometry of compact linear groups
  of cohomogeneity $\leq3$}, {Diff. Geom. and its Appl.} \textbf{4} (1994),
  1--23.

\bibitem[Str96]{Straume2}
E.~Straume, \emph{Compact connected lie transformation groups on spheres with
  low cohomogeneity, {I}}, Memoirs, no. 569, Amer. Math. Soc., 1996.

\bibitem[Tho88]{Th4}
G.~Thorbergsson, \emph{Homogeneous spaces without taut embeddings}, {Duke Math.
  J.} \textbf{57} (1988), 347--355.

\bibitem[TT97]{TTh1}
C.-L. Terng and G.~Thorbergsson, \emph{Taut immersions into complete
  {R}iemannian manifolds}, Tight and Taut Submanifolds (T.~E. Ryan and S.-S.
  Chern, eds.), Math. Sci. Res. Inst. Publ. 32, Cambridge University Press,
  1997, pp.~181--228.

\bibitem[Wol84]{Wolf}
J.A. Wolf, \emph{Spaces of constant curvature}, 5th ed., Publish or Perish,
  Houston, 1984.

\end{thebibliography}

\end{document}